\documentclass[11pt]{amsart}
\usepackage{amsmath,amssymb,amsthm}
\usepackage[margin=1in]{geometry}

\newtheorem{theorem}{Theorem}[section]
\newtheorem{proposition}[theorem]{Proposition}
\newtheorem{lemma}[theorem]{Lemma}
\newtheorem{corollary}[theorem]{Corollary}
\theoremstyle{definition}

\theoremstyle{remark}
\newtheorem{remark}[theorem]{Remark}
\newtheorem{example}[theorem]{Example}
\theoremstyle{plain}
\newtheorem{introtheorem}{Theorem}

\newtheorem{introtheoremprime}{Theorem}

\DeclareMathOperator{\sign}{sign}
\DeclareMathOperator{\cone}{cone}
\DeclareMathOperator{\conv}{conv}
\newcommand{\R}{\mathbb{R}}
\newcommand{\A}{\mathcal A}

\title{Descartes' rule of signs for the topological degree}
\author{Fr\'ed\'eric Bihan}
\address{Laboratoire de Math\'ematiques, Universit\'e Savoie Mont Blanc, 73376 Le Bourget-du-Lac Cedex, France}
\email{frederic.bihan@univ-smb.fr}
\date{\today}

\begin{document}
\begin{abstract}
We consider maps $F\colon\R^n_{>0}\to\R^n$ whose components are real polynomials, possibly with real exponents, with the same support of $N=n+k+1$ monomials, and we study the topological degree of $F$ on the positive orthant. In a previous paper we computed this topological degree in terms of signs of maximal minors, under a half-space condition which holds for generic exponents and coefficients. Here we use these formulas, together with their translation on the Gale dual side, to bound the topological degree. Our main result is a Descartes' rule of signs for the topological degree, which is as simple as the classical one: for two variables, the absolute value of the topological degree is at most half the number of cyclic sign changes, at the vertices of the Newton polygon, of the coefficients of any linear combination of the two equations, and similar statements hold in any dimension and on the Gale dual side. Our second result is a fewnomial bound: for generic exponents and coefficients, we get $|\deg F|<\frac12\min(M(n,N),M(k,N))$, where $M(d,N)$ is the maximal number of facets of a $d$-dimensional polytope with $N$ vertices, and this bound is sharp when $n=2$ or $k=2$. Our third result is a symmetry: using the H-duality of Bihan, Dickenstein and Forsg\aa rd, we show that the set of possible values of the topological degree is the same for the parameters $(n,k)$ and $(k,n)$.
\end{abstract}

\maketitle

\section*{Introduction}

Descartes' rule of signs, which bounds the number of positive roots of a real univariate polynomial by the number of sign changes in the sequence of its coefficients, is the model of the upper bounds for the number of real solutions of polynomial equations. In particular, the number of positive roots is bounded in terms of the number of monomials, independently of the degree. Fewnomial theory, developed by Khovanskii \cite{Kh}, extends this last fact to systems: the number of nondegenerate positive solutions of a system of $n$ polynomial equations in $n$ variables with $N$ monomials is bounded by a function of $n$ and $N$ only. Khovanskii's bound is far from being sharp. It has been improved by several authors (see \cite{BS,BD,BDF,PartI} for instance), but the sharp bounds are known only in very few cases. For instance, the maximal number of positive solutions of a system of two equations in two variables with five monomials is not known, and the best construction, due to El Hilany \cite{EH}, has $7$ positive solutions.

In this paper we consider the same question for a coarser invariant, the topological degree. Let $\A=\{a_1,\dots,a_N\}\subset\R^n$ be a finite set of real exponent vectors which affinely spans $\R^n$, write $N=n+k+1$, and let $C$ be a real $n\times N$ matrix with columns $C_1,\dots,C_N$. We study the topological degree on the positive orthant of the fewnomial map
\[
F=(f_1,\dots,f_n)\colon\R^n_{>0}\to\R^n,\qquad F(x)=\sum_{i=1}^NC_i\,x^{a_i}.
\]
When the positive zeros of $F$ are nondegenerate, $\deg F$ is the sum of the signs of the Jacobian determinant of $F$ at these zeros. Thus $|\deg F|$ is a lower bound for the number of positive zeros, with the same parity. The topological degree can be seen as the robust part of the positive solutions. It does not change when the coefficients vary, as long as no positive solution goes to infinity, while the positive solutions can appear or disappear in pairs with opposite signs. For this reason, the topological degree is a classical tool to prove the existence of positive solutions, see for instance Bihan, Dickenstein and Giaroli \cite{BDG}, and it is used in reaction network theory to count steady states and to detect multistationarity (see \cite{CHW,CFMW}). Degree arguments can only guarantee $|\deg F|$ positive solutions, so an upper bound for $|\deg F|$ limits the lower bounds for the number of positive solutions which can be obtained in this way. It also says how far the number of positive solutions can be from the number which is forced topologically.

In \cite{Degree1} we showed that, under a half-space condition (H) on the columns of $C$ indexed by the faces of the Newton polytope, which is satisfied for generic exponents and coefficients, $\deg F$ is determined by the signs of the maximal minors of $C$ and of the matrix $A$ whose columns are the vectors $(1,a_i)$. One advantage of the topological degree over the number of positive solutions is that it is given by such combinatorial formulas, so that sharp bounds can be obtained in cases where the sharp bounds for the number of positive solutions are not known. Already for $n=1$, a polynomial with $N$ monomials can have $N-1$ positive roots, while its topological degree is $0$ or $\pm1$. For $n=k=2$, we will see that the topological degree is at most $2$ in absolute value, while the system of El Hilany mentioned above has $7$ positive solutions.

Our main tool, besides the formulas of \cite{Degree1}, is Gale duality for polynomial systems, in the sense of Bihan and Sottile \cite{BSGale,BS}. A system with $n$ variables and $N=n+k+1$ monomials is equivalent to a system of $k$ equations in $k$ variables, involving products of powers of affine functions, on a polyhedron. This duality is particularly efficient when $k$ is small, and it was used in \cite{BS} to obtain fewnomial bounds which are polynomial in $n$ for fixed $k$. We translate the formulas of \cite{Degree1} on the Gale dual side, and the two sides together give bounds in terms of both $n$ and $k$.

\subsection*{Our results}
We present three main results: a Descartes' rule of signs for the topological degree, fewnomial bounds for the topological degree, and a symmetry property of the set of possible values of the topological degree.

Let us first recall the classical Descartes' rule of signs. Let $f(x)=c_1x^{a_1}+\dots+c_Nx^{a_N}$ be a real polynomial in one variable, with nonzero coefficients and real exponents $a_1<\dots<a_N$. Then
\[
\#\{\text{positive roots of }f\}\le V(c_1,\dots,c_N),
\]
where the roots are counted with multiplicity and $V(c_1,\dots,c_N)$ is the number of sign changes in the sequence $c_1,\dots,c_N$. Moreover, the difference is even. This bound is remarkably simple: one only has to read the signs of the coefficients in the order of the exponents. Generalizing it to polynomial systems is a long-standing problem. Itenberg and Roy proposed a multivariate Descartes' rule of signs \cite{IR}, which turned out to be false \cite{LW}. Descartes' rules of signs are known for systems supported on circuits, that is, for $k=1$ \cite{BD,BDF}. In \cite{PartI}, we obtained a Descartes' rule of signs for arbitrary fewnomial systems, in the sense that the bounds only depend on sign conditions on the maximal minors of the coefficient and exponent matrices. However, these bounds are given by a recursion over the flats of a hyperplane arrangement, and they are far from the simplicity of the classical rule.

The main result of this paper is that the situation is much better for the topological degree: we obtain a Descartes' rule of signs which is as simple as the classical one. In one variable, the topological degree of $f$ on $\R_{>0}$ only depends on the signs of the coefficients of the extreme monomials, since $\deg f=\frac12(\sign c_N-\sign c_1)$. In two variables, let $v_1,\dots,v_r$ be the vertices of the Newton polygon in cyclic order, and let $g$ be a linear combination of $f_1,f_2$ whose coefficients $g_{v_1},\dots,g_{v_r}$ at these vertices are nonzero. We show (Theorem~\ref{thm:desc-intro}) that
\[
|\deg F|\le\tfrac12V_{\rm cyc}(g_{v_1},\dots,g_{v_r}),
\]
where $V_{\rm cyc}(g_{v_1},\dots,g_{v_r})$ is the number of cyclic sign changes, that is, the number of indices $j\in\{1,\dots,r\}$ such that $g_{v_j}g_{v_{j+1}}<0$, with $v_{r+1}=v_1$. Equivalently, this is the number of edges of the Newton polygon at whose endpoints the coefficients of $g$ have opposite signs. Thus, as in the classical rule, it suffices to read the signs of the coefficients of one polynomial in a suitable order. Moreover, only the vertices of the Newton polygon play a role, as for the extreme monomials in one variable. We find it remarkable that such a simple rule exists, in contrast with the case of the number of positive solutions.

Similar statements hold in any dimension. For $n$ variables, we take $n-1$ linear combinations $g_1,\dots,g_{n-1}$ of the equations, and the role of the sign changes is played by the facets of $P$ (or the simplices of a triangulation of its boundary) on which the truncated system $g_1=\dots=g_{n-1}=0$ has a positive solution. Equivalently, the vectors of coefficients of $g_1,\dots,g_{n-1}$ at the vertices of such a facet positively span $\R^{n-1}$, which for $n=2$ means that the two coefficients have opposite signs. Then $|\deg F|$ is at most half the number of these facets. There is also a dual rule on the Gale side, where the roles of the exponents and of the coefficients are exchanged. There, the linear combinations of the equations are replaced by affine relations among the exponents, and the facets of $P$ by the vertices of the polyhedron of the Gale dual system, which is determined by the coefficients. For $k=2$, the bound is half the number of cyclic sign changes of the values of an affine relation of $\A$, read in the cyclic order of the edges of this polygon. This is a first instance of the symmetry between $n$ and $k$ discussed below.

The key observation is the following. By the formulas of \cite{Degree1}, $\deg F$ is a signed count of the preimages of a direction $u$ by a piecewise linear map from the boundary of the Newton polytope $P=\conv(\A)$ to the sphere. Taking $u=e_n$ and $u=-e_n$ and adding the two formulas, we get that $2\deg F$ is a signed count of the ends of the curve $f_1=\dots=f_{n-1}=0$ (see Proposition~\ref{prop:ends-newton} and the comment after it). On the Gale dual side, the topological degree becomes a signed count of the vertices of the polyhedron of the Gale dual system (Theorem~\ref{thm:B}), and we get a similar statement there (Proposition~\ref{prop:ends-gale}).

Let us now state our Descartes' rule of signs more precisely. Consider $n-1$ linear combinations $g_1,\dots,g_{n-1}$ of $f_1,\dots,f_n$, and a triangulation $T$ of the boundary of $P$ with vertices in $\A$ such that every face of $P$ is a union of simplices of $T$. For an $(n-1)$-simplex $\tau$ of $T$, we denote by $g_j|_\tau$ the sum of the terms of $g_j$ whose exponents are in $\tau$. We say that $\tau$ is a \emph{sign change} of $(g_1,\dots,g_{n-1})$ if the truncated system $g_1|_\tau=\dots=g_{n-1}|_\tau=0$ has a positive solution, and we denote by $V_T(g_1,\dots,g_{n-1})$ the number of these simplices. For $n=2$ and $T$ given by the edges of the Newton polygon, an edge is a sign change of $g$ if and only if the coefficients of $g$ at its two endpoints have opposite signs. We also need a mild genericity condition. For $i\in[N]$, let $\gamma_i\in\R^{n-1}$ be the vector of the coefficients of $x^{a_i}$ in $g_1,\dots,g_{n-1}$. We say that $(g_1,\dots,g_{n-1})$ is \emph{admissible} for $T$ if there is no vector $\lambda\ge0$, supported on at most $n-1$ vertices of $T$, such that
\[
\sum_i\lambda_i\gamma_i=0\qquad\text{and}\qquad\sum_i\lambda_iC_i\neq0 .
\]
This holds for generic linear combinations, and for $n=2$ it only means that the coefficients of $g$ at the vertices of $T$ are nonzero.

\begin{introtheorem}[Descartes' rule of signs for the topological degree]\label{thm:desc-intro}
Assume (H).
\begin{enumerate}
\item For every triangulation $T$ of $\partial P$ as above, and for all linear combinations $g_1,\dots,g_{n-1}$ of $f_1,\dots,f_n$ which are admissible for $T$,
\[
|\deg F|\le\tfrac12V_T(g_1,\dots,g_{n-1}).
\]
\item In particular, if $n=2$ and $v_1,\dots,v_r$ are the vertices of the Newton polygon in cyclic order, then for every linear combination $g=\lambda_1f_1+\lambda_2f_2$ whose coefficients $g_{v_1},\dots,g_{v_r}$ at these vertices are nonzero, we have
\[
|\deg F|\le\tfrac12V_{\rm cyc}(g_{v_1},\dots,g_{v_r}),
\]
where $V_{\rm cyc}$ is the number of cyclic sign changes.
\end{enumerate}
\end{introtheorem}

This is a part of Theorem~\ref{thm:descartes}, which also contains the analogous statements on the Gale dual side, where the role of the simplices of $T$ is played by the vertices of the polyhedron of the Gale dual system, and the role of the linear combinations $g_j$ by affine relations on $\A$. For $n=1$, the triangulation $T$ consists of the two endpoints of $P$, and the first statement gives $|\deg F|\le1$. For $n=2$, if some linear combination of the two equations has coefficients of the same sign at all the vertices of the Newton polygon, then $\deg F=0$. For the map of Example~\ref{ex:run}, the coefficients of $f_1$ at the vertices of the Newton polygon are $3,-3,1,1,-3$. This sequence has $4$ cyclic sign changes, and so $|\deg F|\le2$, which is an equality.

Our second main result is a fewnomial bound for the topological degree. Counting the number of terms in the previous formulas with the Upper Bound Theorem \cite{McM}, we obtain the following bound, which only depends on $n$ and $N$, or on $k$ and $N$.

\begin{introtheorem}\label{thm:C}
Assume that all the maximal minors of $A$ and of $C$ are nonzero. If $n,k\ge2$, then
\[
|\deg F|<\tfrac12\min\big(M(n,N),M(k,N)\big),\qquad M(d,N)=\binom{N-\lceil d/2\rceil}{\lfloor d/2\rfloor}+\binom{N-\lfloor d/2\rfloor-1}{\lceil d/2\rceil-1},
\]
where $M(d,N)$ is the number of facets of a cyclic $d$-dimensional polytope with $N$ vertices. If $n=1$ or $k=1$, then $|\deg F|\le1$.
\end{introtheorem}

Let $\mathcal D(n,k)$ be the set of the values of $|\deg F|$ for the fewnomial maps with $n$ components and $N=n+k+1$ monomials whose exponent and coefficient matrices have nonzero maximal minors. When $n=2$ or $k=2$, we have $\min(M(n,N),M(k,N))=N$, and the bound of Theorem~\ref{thm:C} is sharp.

\begin{introtheoremprime}\label{thm:E}
For all $n,k\ge1$,
\[
\mathcal D(n,1)=\mathcal D(1,n)=\{0,1\},\qquad \mathcal D(n,2)=\mathcal D(2,n)=\Big\{0,1,\dots,\Big\lfloor\frac n2\Big\rfloor+1\Big\}.
\]
\end{introtheoremprime}

In particular, for $n=k=2$ the topological degree is at most $2$ in absolute value, which is much smaller than the $7$ positive zeros of the system of El Hilany \cite{EH}. We do not know the maximal topological degree when $n,k\ge3$. The first open case is $n=k=3$, where Theorem~\ref{thm:C} gives $|\deg F|\le4$.

Our third main result is a symmetry between the dimension $n$ and the codimension $k$. Let $\Xi_{n,k}$ be the maximal number of nondegenerate positive solutions of a system of $n$ real polynomial equations in $n$ variables with at most $n+k+1$ monomials. Bihan, Santos and Spaenlehauer \cite{BSS} asked whether $\Xi_{n,k}=\Xi_{k,n}$. This symmetry holds for all the known exact values and lower bounds, but it is open in general. A similar question was asked by Bihan, Dickenstein and Forsg\aa rd \cite{BDF24} for the maximal multiplicity of an isolated solution in the complex torus of such a system. To study it, they introduced the H-duality, which associates with a system with $n$ variables and codimension $k$ a system with $k$ variables and codimension $n$, and they showed that the multiplicity is preserved under convenience and nondegeneracy assumptions. They answered their question positively when $k=1$.

Here we use the H-duality to prove the analogous symmetry for the topological degree. Note first that Theorems~\ref{thm:C} and~\ref{thm:E} are symmetric in $n$ and $k$, and this is not a coincidence. The H-dual of $F$ is a fewnomial map $H$ with $k$ components in $k$ variables and $N$ monomials. The two maps have in general different numbers of positive zeros, but the H-duality exchanges the Newton polytope formula for $F$ with the Gale dual formula for $H$, and we get that $\deg H=\pm\deg F$ (Theorem~\ref{thm:h-duality}). As a consequence, we obtain the following result.

\begin{introtheorem}\label{thm:D}
For all $n,k\ge1$, we have $\mathcal D(n,k)=\mathcal D(k,n)$.
\end{introtheorem}

In other words, the answer to the question of \cite{BSS} is positive if we replace the number of positive solutions by the topological degree. The question for $\Xi_{n,k}$ itself remains open.

\subsection*{Previous work}
In \cite{PartI}, upper bounds for the number of positive solutions are obtained by a Khovanskii--Rolle argument on the Gale side. There, the number of noncompact components of a curve bounds $|\deg F|$ (see the discussion after Proposition~\ref{prop:ends-gale}). Bihan, Dickenstein and Giaroli \cite{BDG} combine degree theory and Gale duality to guarantee the existence of positive solutions.

\subsection*{Organization}
Section~\ref{sec:background} recalls the definitions and the results of \cite{Degree1}. Section~\ref{sec:gale} proves Theorem~\ref{thm:B} and treats the cases $k\le2$. Section~\ref{sec:ends} writes the topological degree as a half signed count, and Section~\ref{sec:descartes} deduces our Descartes' rule of signs. Section~\ref{sec:bound} proves Theorems~\ref{thm:C} and~\ref{thm:E}, and Section~\ref{sec:h} the H-duality and Theorem~\ref{thm:D}.

\subsection*{Acknowledgements and use of AI}
This paper grew out of \cite{PartI}. The author used the AI assistant Claude Opus~5.5 (Anthropic) for help with several technical proofs, for checking the examples, and for help with the writing. The author checked all the mathematical content and takes full responsibility for it.

\subsection*{Data availability}
Data sharing is not applicable to this article, as no datasets were generated or analysed during the current study.

\subsection*{Conflict of interest}
The author has no conflicts of interest to declare that are relevant to the content of this article.

\section{Background}\label{sec:background}

In this section we recall the notation and the results of \cite{Degree1} that we use. We first recall the half-space condition (H), under which the topological degree of $F$ is well defined and does not change when the coefficients vary. We then recall the two formulas of \cite{Degree1} for the topological degree. The Newton polytope formula expresses $\deg F$ as the degree of a piecewise linear map from the boundary of the Newton polytope to the sphere, and it is the starting point of Section~\ref{sec:ends}. The Viro formula expresses $\deg F$ as a signed count of positively decorated simplices of a regular triangulation, and it is translated on the Gale dual side in Section~\ref{sec:gale}.

\subsection{The topological degree}
For $\tau\subset[N]$, $A^\tau$ and $C^\tau$ denote the submatrices formed by the columns indexed by $\tau$, in increasing order. We identify a face $G$ of the Newton polytope $P$ with the set of indices $i$ such that $a_i\in G$, and we write $F_G=\sum_{i\in G}C_ix^{a_i}$. We say that (U$_A$), respectively (U$_C$), holds if all the maximal minors of $A$, respectively of $C$, are nonzero. We consider the \emph{half-space condition}
\begin{itemize}
\item[(H)] for every proper face $G$ of $P$, the columns $C_i$, $i\in G$, lie in an open half-space of $\R^n$.
\end{itemize}
By Gordan's theorem, (H) means that $\ker C^G$ contains no nonzero vector with nonnegative coordinates. It holds under (U$_A$) and (U$_C$), and it is preserved by the rescalings $C\mapsto C\,\mathrm{diag}(\lambda)$, $\lambda\in\R^N_{>0}$. Under (H), no truncation $F_G$ to a proper face has a positive zero, the positive zeros of $F$ form a compact set, and $\deg F=\deg(F,\R^n_{>0},0)$ is the Brouwer degree of $F$ on $\R^n_{>0}$. It does not change along continuous families of coefficient matrices satisfying (H). If the positive zeros are nondegenerate, their number is at least $|\deg F|$ and congruent to $\deg F$ modulo $2$ \cite[Section~1]{Degree1}.

\subsection{The Newton polytope formula}\label{sec:onetwo}
For a triangulation $T$ of $\partial P$ with vertices in $\A$, let $\Pi_T$ be the continuous map which is affine on every simplex of $T$ and sends each vertex $a_s$ of $T$ to $C_s$. For an $(n-1)$-simplex $\tau$ of $T$, let $i_\tau$ be an index such that $a_{i_\tau}$ is off the facet containing $\tau$. The sign of $\det A^{\tau\cup\{i_\tau\}}$, where the column $i_\tau$ is placed last, does not depend on this choice.

\begin{lemma}[{\cite[Lemmas~3.1 and~3.2]{Degree1}}]\label{lem:homotopy}\label{lem:PL}
Assume (H), and let $T$ be a triangulation of $\partial P$ with vertices in $\A$ such that every face of $P$ is a union of simplices of $T$. Then $\Pi_T$ does not vanish on $\partial P$, and $\deg F$ is the degree of $g=\Pi_T/|\Pi_T|\colon\partial P\to S^{n-1}$, where $\partial P$ is oriented as the boundary of $P$. Let $u\in\R^n$ be a vector which does not lie in any cone spanned by $n-1$ columns $C_s$, $s$ a vertex of $T$. An $(n-1)$-simplex $\tau$ of $T$ meets $g^{-1}(u)$ if and only if $u\in\cone(C_s,\ s\in\tau)$, and then in exactly one point, which lies in the relative interior of $\tau$, and the local degree of $g$ at this point is $(-1)^n\sign\det C^\tau\cdot\sign\det A^{\tau\cup\{i_\tau\}}$.
\end{lemma}

\begin{theorem}[{\cite[Theorem~A$'$]{Degree1}}]\label{thm:Ap}
Assume (H), and let $T$ and $u$ be as in Lemma~\ref{lem:PL}. Then
\[
\deg F=(-1)^n\sum_{\substack{\tau\in T,\ \dim\tau=n-1\\ u\in\cone(C_s,\ s\in\tau)}}\sign\det C^\tau\cdot\sign\det A^{\tau\cup\{i_\tau\}} .
\]
\end{theorem}

Under (U$_A$), every facet $\Gamma$ of $P$ is a simplex which contains exactly $n$ points of $\A$, and we can take for $T$ the boundary complex of $P$. We write $\varepsilon_\Gamma=\sign\det A^{\Gamma\cup\{i_\Gamma\}}$.

\begin{theorem}[{\cite[Theorem~A]{Degree1}}]\label{thm:A}
Assume (U$_A$) and (U$_C$), and let $u\in\R^n$ be a vector which does not lie in any cone spanned by $n-1$ columns of $C$. Then
\[
\deg F=(-1)^n\sum_{\substack{\Gamma\ \text{facet of}\ P\\ u\in\cone(C_s,\ s\in\Gamma)}}\sign\det C^\Gamma\cdot\varepsilon_\Gamma .
\]
\end{theorem}

\begin{theorem}[{\cite[Theorem~3.5]{Degree1}}]\label{cor:om}
Assume (H). Then $\deg F$ only depends on the oriented matroids of the columns of $A$ and of $C$ indexed by the vertices of $P$.
\end{theorem}

For $n=1$, with $a_1<\dots<a_N$ and $F=\sum_ic_ix^{a_i}$, Theorem~\ref{thm:A} gives $\deg F=\frac12(\sign c_N-\sign c_1)$ \cite[Section~3.1]{Degree1}. For $n=2$ we have the following.

\begin{corollary}[{\cite[Corollary~3.6]{Degree1}}]\label{cor:n2}
Let $n=2$, assume (H), and let $j_1,\dots,j_m$ be the indices of the vertices of the Newton polygon $P$, in counterclockwise order. Then the closed polygonal line with successive vertices $C_{j_1},C_{j_2},\dots,C_{j_m}$, that is, the union of the segments $[C_{j_r},C_{j_{r+1}}]$ with $j_{m+1}=j_1$, does not pass through the origin, and $\deg F$ is its winding number around the origin.
\end{corollary}

\subsection{The Viro formula}
For $h\in\R^N$, a set $\sigma\subset[N]$ with $|\sigma|=n+1$ is a simplex of the regular subdivision $T_h$ if there is an affine function $\lambda_\sigma$ on $\R^n$ such that $\lambda_\sigma(a_i)=h_i$ for $i\in\sigma$ and $\lambda_\sigma(a_i)<h_i$ for $i\notin\sigma$, which means that $\{(a_i,h_i)\}_{i\in\sigma}$ spans a lower facet of the convex hull of the lifted points. We say that $h$ is \emph{generic} if every lower facet contains exactly $n+1$ lifted points; then $T_h$ is a triangulation of $P$. A real $n\times(n+1)$ matrix $M$ is \emph{positively spanning} if the numbers $(-1)^j\operatorname{minor}(M,j)$, $j=1,\dots,n+1$, are nonzero and have the same sign, where $\operatorname{minor}(M,j)$ is the determinant of the matrix obtained by deleting the $j$-th column of $M$. When $M$ has rank $n$, these numbers are the coordinates of a vector spanning $\ker M$. Following \cite{BSS}, a simplex $\sigma\in T_h$ is \emph{positively decorated} by $C$ if $C^\sigma$ is positively spanning, and its sign is $s(\sigma)=\sign\det A^\sigma\cdot\sign\det\widehat C^\sigma$, where $\widehat C^\sigma$ is the matrix with columns $(1,C_l)$, $l\in\sigma$, in increasing order.

\begin{theorem}[{\cite[Theorem~2.2]{Degree1}}]\label{thm:viro}
Assume (U$_C$) and (H), and let $h\in\R^N$ be generic.
\begin{enumerate}
\item For $t>0$ small enough, the positive zeros of the Viro map $F^{h,t}(x)=\sum_it^{h_i}C_ix^{a_i}$ are nondegenerate, and there is exactly one of them for each positively decorated simplex $\sigma\in T_h$, and no other. The sign of the Jacobian determinant at the zero attached to $\sigma$ is $s(\sigma)$.
\item Consequently,
\[
\deg F=\sum_{\substack{\sigma\in T_h\\ \sigma\ \text{positively decorated}}}s(\sigma).
\]
\end{enumerate}
\end{theorem}

\begin{example}[A running example]\label{ex:run}
Let
\[
F(x,y)=\big(3x-3x^3+x^4y^2+x^2y^3-3y^2,\ 2x^3-3x^4y^2+3x^2y^3-2y^2\big),
\]
with $\A=\{(1,0),(3,0),(4,2),(2,3),(0,2)\}$, the vertices of a pentagon in counterclockwise order, and columns $C_1=(3,0)$, $C_2=(-3,2)$, $C_3=(1,-3)$, $C_4=(1,3)$, $C_5=(-3,-2)$. Here $n=k=2$, and (U$_A$) and (U$_C$) hold. The closed polygonal line with successive vertices $C_1,\dots,C_5$ turns twice around the origin, so $\deg F=2$. The heights $h=(-1,0,0,0,0)$ induce the fan $T_h=\{\{1,2,3\},\{1,3,4\},\{1,4,5\}\}$ from $a_1$, in which $\{1,2,3\}$ and $\{1,4,5\}$ are positively decorated, both with sign $+1$ \cite[Examples~2.4 and~3.7]{Degree1}.
\end{example}

\section{The Gale dual formula}\label{sec:gale}

In this section we translate the formulas of \cite{Degree1} on the Gale dual side. Recall that Gale duality \cite{BSGale,BS} transforms the positive solutions of $F$ into the solutions of a system of $k$ equations in $k$ variables in a polyhedron $\Delta\subset\R^k$, which is determined by the coefficient matrix $C$ (see~\eqref{eq:gale-system}). The main result of this section is the Gale dual formula (Theorem~\ref{thm:B}), which expresses $\deg F$ as a signed count of some vertices of $\Delta$, selected by a vector $e$. This formula is the Viro formula written on the Gale side: the simplices of a regular triangulation of $\A$ correspond to the cones spanned by rows of a Gale dual matrix of $A$ which contain $e$, and the positively decorated simplices correspond to the vertices of $\Delta$ (Lemma~\ref{lem:translation}). In Section~\ref{sec:gale-matrices} we introduce the Gale dual matrices and the polyhedron $\Delta$, and we state Theorem~\ref{thm:B}. We prove it in Section~\ref{sec:gale-viro}, and in Section~\ref{sec:codim} we make it explicit when $k=1$ or $k=2$. The Gale dual formula is used in Section~\ref{sec:ends} to write the topological degree as a half signed count on the Gale side, and in Section~\ref{sec:h} to prove the symmetry between $n$ and $k$.

\subsection{Gale dual matrices}\label{sec:gale-matrices}
We first fix the notation for the Gale dual matrices, then we introduce the polyhedron $\Delta$ and its vertices, and we state the Gale dual formula.

Let $B$ be an $N\times k$ matrix whose columns form a basis of $\ker A$, and let $D$ be an $N\times(k+1)$ matrix whose columns form a basis of $\ker C$. We say that $B$ and $D$ are \emph{Gale dual matrices} of $A$ and $C$. Since the first row of $A$ is $(1,\dots,1)$, the rows of $B$ sum to zero: $\sum_iB_i=0$. We denote by $B_i\in\R^k$ and $D_i\in\R^{k+1}$ the $i$-th rows, and for $S\subset[N]$ we denote by $B_S$ and $D_S$ the submatrices with rows in $S$, in increasing order. Under (U$_C$), $D_N\neq0$, and we normalize $D$ so that $D_N=(0,\dots,0,1)$.

We will use the following classical relation between complementary minors \cite[Section~3.4]{BLSWZ}. If $M$ is an $r\times N$ matrix of rank $r$ and the columns of the $N\times(N-r)$ matrix $G$ form a basis of $\ker M$, then there is a constant $c\neq0$ such that
\begin{equation}\label{eq:complementary}
\det M^\tau=c\,(-1)^{\sum_{j\in\tau}j}\,\det G_{[N]\setminus\tau}\qquad\text{for all }\tau\subset[N],\ |\tau|=r .
\end{equation}
In particular, (U$_A$) is equivalent to the condition that any $k$ rows of $B$ are linearly independent, and (U$_C$) to the condition that any $k+1$ rows of $D$ are linearly independent.

We put $P_i(v)=\langle D_i,v\rangle$ for $v\in\R^{k+1}$ and $p_i(y)=P_i(y,1)$ for $y\in\R^k$, so that $p_N=1$, and we consider
\[
\Delta=\{y\in\R^k:\ p_i(y)>0,\ i\in[N]\},\qquad K=\{v\in\R^{k+1}:\ P_i(v)\ge0,\ i\in[N]\}.
\]
The vectors of $\ker C$ are the vectors $Dv$. It follows that $\Delta\neq\emptyset$ if and only if $\ker C$ contains a positive vector. When $\Delta\neq\emptyset$, $K$ is a pointed cone with nonempty interior, and $\overline\Delta$ is identified with the image of $K\setminus\{0\}$ in the sphere $S^k$. The polyhedron $\Delta$ carries the Gale dual system of Bihan and Sottile \cite{BSGale}. Namely, the map $x\mapsto y$ defined by $p_i(y)=x^{a_i-a_N}$ is a bijection from the positive zeros of $F$ onto the zeros in $\Delta$ of
\begin{equation}\label{eq:gale-system}
\psi=(\psi_1,\dots,\psi_k),\qquad\psi_j=\sum_iB_{ij}\log p_i .
\end{equation}
Indeed, $F(x)=0$ with $x>0$ means that the positive vector $(x^{a_i-a_N})_i$, whose last coordinate is $1$, lies in $\ker C$, so it is equal to $D(y,1)$ for a unique $y\in\Delta$. Moreover, a positive vector of the form $D(y,1)$ is of the form $(x^{a_i-a_N})_i$ if and only if its logarithm lies in the row space of $A$, which means that $\psi(y)=0$.

For $S\subset[N]$, let $K_S=K\cap\{P_s=0,\ s\in S\}$. We say that $S$ is a \emph{face} of $\Delta$ if $K_S\neq\{0\}$ and $P_i>0$ on the relative interior of $K_S$ for $i\notin S$. We denote by $V(\Delta)$ the set of the faces of size $k$. They correspond to the vertices of the closure $\overline\Delta$ of $\Delta$ in the projective space, or equivalently to the extreme rays of $K$, including the vertices at infinity when $\Delta$ is unbounded.

\begin{proposition}\label{prop:vertices}
Assume (U$_C$) and $\Delta\neq\emptyset$.
\begin{enumerate}
\item Every face $S$ of $\Delta$ has $|S|\le k$, and every subset of a face is a face. The sets $S\in V(\Delta)$ correspond to the extreme rays of $K$.
\item A set $S$ with $|S|=k$ belongs to $V(\Delta)$ if and only if the numbers $\det(D_S;D_i)$, $i\notin S$, all have the same sign.
\end{enumerate}
\end{proposition}

\begin{proof}
By (U$_C$), any $k+1$ of the linear forms $P_i$ are linearly independent. Then a nonzero face of $K$ lies on at most $k$ of the hyperplanes $\{P_i=0\}$. For $|S|=k$, the line $\{P_s=0,\ s\in S\}$ is spanned by the cofactor vector $v_S$ of $D_S$, which satisfies $\langle D_i,v_S\rangle=\det(D_S;D_i)$. This line meets $K$ in an extreme ray if and only if $\pm v_S\in K$, which means that the nonzero numbers $\det(D_S;D_i)$, $i\notin S$, have the same sign. Finally, every nonzero face of $K$ is spanned by its extreme rays, and the rest follows from the simplicity of the arrangement.
\end{proof}

Let $\lambda=\sum_iD_i$, so that $\lambda(v)=\sum_iP_i(v)>0$ on $K\setminus\{0\}$, and let $Q=\{v\in K:\ \lambda(v)=1\}$. This is a compact simple polytope of dimension $k$, which is a compact model of the closure $\overline\Delta$ of $\Delta$ in the projective space: a point $y\in\Delta$ corresponds to the point $(y,1)/\lambda(y,1)$ of the interior of $Q$, and the points at infinity of $\overline\Delta$ correspond to the points of $Q$ with last coordinate $0$. Its vertices are the points $Q\cap K_S$, $S\in V(\Delta)$, and its facets lie on the hyperplanes $\{P_i=0\}$. In particular, $Q$ has at most $N$ facets.

For $S\in V(\Delta)$, choose an index $i_S\notin S$ and put
\[
\omega(S)=\sign\det B_S\cdot\sign\det(D_S;D_{i_S}),
\]
where $B_S$ is formed by the rows of $B$ indexed by $S$, and $(D_S;D_{i_S})$ by the rows of $D$ indexed by $S$ followed by the row $i_S$. By Proposition~\ref{prop:vertices}, this sign does not depend on the choice of $i_S$.

\begin{theorem}[Gale dual formula]\label{thm:B}\label{thm:Bp}
Assume (U$_C$) and (H), and let $e\in\R^k$ be a vector which does not belong to any cone spanned by $k-1$ of the rows $-B_1,\dots,-B_N$. If $\Delta=\emptyset$, then $\deg F=0$. Otherwise,
\[
\deg F=\pm\sum_{\substack{S\in V(\Delta)\\ e\in\cone(-B_s,\ s\in S)}}\omega(S),
\]
where the global sign only depends on the choice of $B$ and $D$.
\end{theorem}

We denote by $\deg(\Delta)$ the sum in Theorem~\ref{thm:B}, and we put $\deg(\Delta)=0$ if $\Delta=\emptyset$. By Theorem~\ref{thm:B}, $\deg(\Delta)=\pm\deg F$, so it does not depend on $e$.

For $k=1$ this gives the parity part of the Descartes' rule of signs for circuits \cite{BD}, and for $k=2$ the topological degree is the winding number of the closed polygonal line whose successive vertices are the rows of $B$ indexed by the edges of $\overline\Delta$ (Proposition~\ref{prop:k2}).

\subsection{Gale translation of the Viro formula}\label{sec:gale-viro}
We now prove Theorem~\ref{thm:B} by translating the Viro formula (Theorem~\ref{thm:viro}). The key point is Lemma~\ref{lem:translation}. A lifting $h\in\R^N$ of $\A$ corresponds to the vector $e=-B^Th\in\R^k$, the simplices of the regular triangulation $T_h$ correspond to the cones spanned by $k$ rows of $-B$ which contain $e$, and the positively decorated simplices correspond to the vertices of $\Delta$.

\begin{lemma}\label{lem:translation}
Assume (U$_C$). Let $h\in\R^N$, put $e=-B^Th$, and let $\sigma\subset[N]$ with $|\sigma|=n+1$ and $S=[N]\setminus\sigma$, so that $|S|=k$.
\begin{enumerate}
\item If $\det B_S\neq0$, then $\sigma\in T_h$ if and only if $e$ lies in the interior of the simplicial cone $\cone(-B_s,\ s\in S)$. If $e$ does not lie in a cone spanned by $k-1$ of the vectors $-B_i$, then $h$ is generic, and the simplices of $T_h$ are the sets $\sigma$ such that $e\in\cone(-B_s,\ s\in S)$; for these, $\det B_S\neq0$. Under (U$_A$), $h$ is generic if and only if $e$ does not lie in a cone spanned by $k-1$ of the vectors $-B_i$.
\item $\sigma$ is positively decorated if and only if $\Delta\neq\emptyset$ and $S\in V(\Delta)$.
\item There is a sign $\epsilon\in\{\pm1\}$, depending only on the choices of $B$ and $D$, such that $\omega(S)=\epsilon\,s(\sigma)$ for every positively decorated $\sigma\in T_h$.
\end{enumerate}
\end{lemma}

\begin{proof}
1. By~\eqref{eq:complementary}, $\det B_S\neq0$ if and only if $\det A^\sigma\neq0$, which means that the points $a_l$, $l\in\sigma$, are affinely independent. Let $\sigma'\subset[N]$ be any subset, with complement $S'$. The points $(a_i,h_i)$, $i\in\sigma'$, span a lower face of the lifted configuration which contains no other lifted point if and only if $h=A^Tc+r$ for some $c\in\R^{n+1}$ and $r\in\R^N$ with $r_{\sigma'}=0$ and $r_{S'}>0$. In this case $B^Th=B^Tr=\sum_{s\in S'}r_sB_s$, so $e=\sum_{s\in S'}r_s(-B_s)$. Conversely, assume that $e=\sum_{s\in S'}r_s(-B_s)$ with $r_s>0$, and extend $r$ by $0$ on $\sigma'$. Then $B^T(h-r)=0$, so $h-r$ lies in the row space of $A$, which is the orthogonal complement of $\ker A$. Taking $\sigma'=\sigma$ with $\det A^\sigma\neq0$, we get the first assertion.

Assume now that $e$ does not lie in a cone spanned by $k-1$ of the vectors $-B_i$. Then a lower facet cannot contain more than $n+1$ lifted points, since its complement $S'$ would have fewer than $k$ elements. So the lower facets are simplices and $h$ is generic. If $\det B_S=0$, then by Carath\'eodory's theorem $\cone(-B_s,\ s\in S)$ is covered by cones spanned by $k-1$ of the vectors $-B_s$, so $e$ does not lie in it. If $\det B_S\neq0$ and $e$ lies in this cone, then it lies in its interior. Finally, under (U$_A$) all the cones $\cone(-B_s,\ s\in S)$ are simplicial. The converse then follows, since $e$ lies in a cone spanned by fewer than $k$ vectors $-B_i$ if and only if it lies in the boundary of some $\cone(-B_s,\ s\in S)$.

2. By (U$_C$), the kernel of $C^\sigma$ is a line. Then $\sigma$ is positively decorated if and only if $\ker C$ contains a vector $m=Dv$ with $m_\sigma>0$ and $m_S=0$, or equivalently such that $P_s(v)=0$ for $s\in S$ and $P_i(v)>0$ for $i\in\sigma$. Moving $v$ slightly in a direction where the $k$ independent forms $P_s$, $s\in S$, are positive, we get an interior point of $K$, so $\Delta\neq\emptyset$. Moreover, $v$ spans an extreme ray of $K$ with $P_i>0$ off $S$, and so $S\in V(\Delta)$. The converse is Proposition~\ref{prop:vertices}.

3. Write $\sigma=\{l_1<\dots<l_{n+1}\}$, let $i\in\sigma$, and let $p$ be the position of $i$ in $\sigma$. Expanding $\det\widehat C^\sigma$ along its first row gives $\sum_{j}(-1)^{1+j}\det C^{\sigma\setminus\{l_j\}}$. Up to a common sign, the terms of this sum are the coordinates of a vector spanning $\ker C^\sigma$, so they have the same sign when $\sigma$ is positively decorated. It follows that $\sign\det\widehat C^\sigma=(-1)^{1+p}\sign\det C^{\sigma\setminus\{i\}}$. By~\eqref{eq:complementary} applied to $(A,B)$ with $\tau=\sigma$, and to $(C,D)$ with $\tau=\sigma\setminus\{i\}$, we have
\[
\sign\det B_S=\pm(-1)^{\sum_{j\in\sigma}j}\sign\det A^\sigma,\qquad \sign\det D_{S\cup\{i\}}=\pm(-1)^{\sum_{j\in\sigma}j-i}\sign\det C^{\sigma\setminus\{i\}},
\]
where the signs $\pm$ only depend on $B$ and $D$. Moreover, $\det(D_S;D_i)=(-1)^{\#\{s\in S:\ s>i\}}\det D_{S\cup\{i\}}$. Multiplying, we get
\[
\omega(S)=\pm(-1)^{i+\#\{s\in S:\ s>i\}+1+p}\,s(\sigma).
\]
Since $i=p+\#\{s\in S:\ s<i\}$, the exponent is congruent to $k+1$ modulo $2$. We get $\omega(S)=\pm(-1)^{k+1}s(\sigma)$, with a sign which does not depend on $\sigma$.
\end{proof}

\begin{proof}[Proof of Theorem~\ref{thm:B}]
We assume (U$_C$) and (H). Since $B^T$ is onto, the generic vector $e$ is of the form $-B^Th$, and $h$ is then generic by Lemma~\ref{lem:translation}(1). If $\Delta=\emptyset$, then no simplex is positively decorated by Lemma~\ref{lem:translation}(2), and $\deg F=0$ by the Viro formula (Theorem~\ref{thm:viro}(2)). Otherwise, by Theorem~\ref{thm:viro}(2) and Lemma~\ref{lem:translation}, we get
\[
\deg F=\sum_{\substack{\sigma\in T_h\\ \text{pos.\ decorated}}}s(\sigma)=\epsilon\sum_{\substack{S\in V(\Delta)\\ e\in\cone(-B_s,\ s\in S)}}\omega(S)=\epsilon\deg(\Delta).\qedhere
\]
\end{proof}

Theorem~\ref{thm:B} is thus the Viro formula written on the Gale side. The simplices of $T_h$ are the complements of the sets $S$ such that $e\in\cone(-B_s,\ s\in S)$, and the positively decorated ones correspond to the vertices of $\overline\Delta$.
\subsection{Codimension one or two}\label{sec:codim}
In this subsection we make Theorem~\ref{thm:B} explicit for small $k$. For $k=1$ we recover the parity part of the Descartes' rule of signs for circuits, and for $k=2$ the topological degree is the winding number of a closed polygonal line whose vertices are rows of $B$ (Proposition~\ref{prop:k2}). This is the Gale dual analogue of the case $n=2$ (Corollary~\ref{cor:n2}).

For $k=1$, $\Delta$ is an interval, and $V(\Delta)$ consists of its two endpoints $\{i_0\}$ and $\{i_1\}$. One of them has index $N$ if $\Delta$ is unbounded, as one sees in the model $Q$. The two signs $\sign\det(D_{i_0};D_i)$ and $\sign\det(D_{i_1};D_i)$ are opposite. Then Theorem~\ref{thm:B} with $e=1$ gives $\deg F=\pm1$ if $B_{i_0}B_{i_1}<0$, and $\deg F=0$ otherwise. This is the parity part of the Descartes' rule of signs for circuits \cite{BD}.

We say that a set $S\subset[N]$ is a \emph{coface} if $S\neq\emptyset$, $S\neq[N]$, and the rows $B_s$, $s\in S$, have a positive linear dependence: $\sum_{s\in S}\lambda_sB_s=0$ for some $\lambda\in\R^S_{>0}$.

\begin{lemma}\label{lem:gale-N}
Let $\A$ and $C$ be arbitrary.
\begin{enumerate}
\item A set $G\subset[N]$ is the set of indices of the points of $\A$ on a nonempty proper face of $P$ if and only if $[N]\setminus G$ is a coface.
\item (H) holds if and only if $K_S=\{0\}$ for every coface $S$; when $\Delta\neq\emptyset$, this means that the closure of $\Delta$ in the projective space does not meet the linear subspace $\{p_s=0,\ s\in S\}$.
\item Under (U$_A$), every coface has at least $k+1$ elements; under (U$_C$), $K_S=\{0\}$ whenever $|S|\ge k+1$.
\end{enumerate}
\end{lemma}

\begin{proof}
1. The vectors $\lambda\in\R^N$ such that $\sum_i\lambda_iB_i=0$ form the row space of $A$, and so they are the vectors $(\ell(a_i))_i$ for the affine functions $\ell$ on $\R^n$. A nonempty proper face is the zero set on $P$ of such a function $\ell$ which is nonnegative on $\A$ and does not vanish identically on $\A$. Its complement in $\A$ is the support of $(\ell(a_i))_i$.

2. By Gordan's theorem, (H) means that $\ker C^G$ contains no nonzero nonnegative vector for every proper face $G$. Extended by $0$ on $S=[N]\setminus G$, such a vector is a vector $Dv$ with $v\in K_S\setminus\{0\}$. We conclude by~1.

3. Under (U$_A$), a proper face contains at most $n$ points of $\A$. Under (U$_C$), any $k+1$ rows of $D$ are linearly independent.
\end{proof}

Lemma~\ref{lem:gale-N}(3) shows again that (U$_A$) and (U$_C$) imply (H).

\begin{proposition}\label{prop:k2}
Let $k=2$, assume (U$_C$) and $\Delta\neq\emptyset$, and let $i_1,\dots,i_m$ be the indices of the edges of the polygon $Q$ (the compact model of $\overline\Delta$ defined in Section~\ref{sec:gale-matrices}), in counterclockwise order. Then (H) holds if and only if the closed polygonal line with successive vertices $B_{i_1},B_{i_2},\dots,B_{i_m}$ does not pass through the origin, and in this case $\deg F$ is, up to sign, the winding number of this polygonal line around the origin.
\end{proposition}

\begin{proof}
Put $i_{m+1}=i_1$. The vertices of $Q$ are the pairs $S_r=\{i_r,i_{r+1}\}$ of consecutive edges.

\emph{The condition (H).} By Lemma~\ref{lem:gale-N}(2), (H) means that $K_S=\{0\}$ for every coface $S$. By Lemma~\ref{lem:gale-N}(3), this only concerns the cofaces with $|S|\le2$, namely the sets $\{i\}$ with $B_i=0$ and the sets $\{a,b\}$ with $\lambda B_a+\mu B_b=0$ for some $\lambda,\mu>0$. By (U$_C$), at most two of the forms $P_i$ vanish at a nonzero point of $K$, and exactly two of them vanish on an extreme ray. It follows that $K_{\{i\}}\neq\{0\}$ if and only if $i$ is an edge of $Q$, since the nonzero face $K_{\{i\}}$ of $K$ contains an extreme ray, on which $P_i$ is one of the two vanishing forms. Likewise, for $|S|=2$, $K_S\neq\{0\}$ if and only if $S$ is one of the vertices $S_r$. As a consequence, (H) holds if and only if no row $B_{i_r}$ vanishes and no two consecutive rows $B_{i_r}$, $B_{i_{r+1}}$ have a positive linear dependence, or equivalently if no segment $[B_{i_r},B_{i_{r+1}}]$ contains the origin.

\emph{The topological degree.} Assume (H). Let $W$ be a $3\times2$ matrix whose columns form a positive basis of $\ker\lambda$, and let $v_0\in Q$. For $a,b\in[N]$, the matrix $(D_a;D_b;\lambda)(W,v_0)$ has rows $(\nu_a,P_a(v_0))$, $(\nu_b,P_b(v_0))$ and $(0,0,1)$, where $\nu_a=D_aW$ is the gradient of $P_a$ restricted to the plane of $Q$. We get $\det(D_a;D_b;\lambda)\det(W,v_0)=\det(\nu_a,\nu_b)$. The vectors $\nu_{i_r}$ are inward normals of the consecutive edges of the convex polygon $Q$, so $\det(\nu_{i_r},\nu_{i_{r+1}})>0$ for all $r$. Moreover, $\det(D_{S_r};\lambda)=\sum_{i\notin S_r}\det(D_{S_r};D_i)$, and all these terms have the same sign. So $\sign\det(D_{S_r};D_i)$ is the same for all $r$, and $\omega(S_r)=\pm\sign\det(B_{i_r},B_{i_{r+1}})$ with a sign which does not depend on $r$. (In $\omega(S_r)$ the rows of $B_{S_r}$ and $D_{S_r}$ are taken in increasing order; taking both in the order $(i_r,i_{r+1})$ changes both signs, so it does not change $\omega(S_r)$. Likewise, the orientation of $\ker\lambda$ only affects a global sign.) Let $e$ be as in Theorem~\ref{thm:B}; here this means that $-e$ does not lie on any ray $\R_{\ge0}B_i$. If $\det(B_{i_r},B_{i_{r+1}})=0$, then by (H) the two rows point in the same direction. Then $e\notin\cone(-B_{i_r},-B_{i_{r+1}})$, and the segment $[B_{i_r},B_{i_{r+1}}]$ does not meet the ray $\R_{>0}(-e)$. Otherwise, $e\in\cone(-B_{i_r},-B_{i_{r+1}})$ if and only if this segment crosses the ray $\R_{>0}(-e)$. By Theorem~\ref{thm:Bp}, we get that $\deg F$ is, up to sign, the number of segments $[B_{i_r},B_{i_{r+1}}]$ which cross the ray $\R_{>0}(-e)$, counted with the sign of $\det(B_{i_r},B_{i_{r+1}})$. This is the winding number.
\end{proof}

Under (U$_A$), any two rows of $B$ are linearly independent. So (H) holds, and the closed polygonal line with successive vertices $B_{i_1},\dots,B_{i_m}$ never passes through the origin. Proposition~\ref{prop:k2} is the Gale dual counterpart of Corollary~\ref{cor:n2}. Note however the following difference. On the Newton side, (H) implies that the closed polygonal line through the columns $C_i$ indexed by the points of $\A$ on $\partial P$ does not pass through the origin, but it is a stronger condition, since it asks that the columns of each edge lie in an open half-plane. On the Gale side, (H) is exactly the condition that the winding number is defined.

For instance, if the edges of a triangle $Q$ correspond to the rows $(1,0)$, $(0,1)$, $(-1,-1)$, then the triangle with these vertices contains the origin, and $|\deg F|=1$. If the last row is $(1,1)$ instead, then $\deg F=0$.

\begin{example}[Example~\ref{ex:run}, continued: the Gale side]\label{ex:run-gale}
Here $k=2$. The vectors $(-5,7,-6,4,0)$ and $(-2,2,-1,0,1)$ form a basis of $\ker A$, and the vectors $(7,9,6,0,0)$, $(-11,-9,0,6,0)$, $(2,1,0,0,1)$ form a basis of $\ker C$. We may thus take
\[
B=\begin{pmatrix}-5&-2\\7&2\\-6&-1\\4&0\\0&1\end{pmatrix},\qquad
D=\begin{pmatrix}7&-11&2\\9&-9&1\\6&0&0\\0&6&0\\0&0&1\end{pmatrix}.
\]
Then $p_1=7y_1-11y_2+2$, $p_2=9y_1-9y_2+1$, $p_3=6y_1$, $p_4=6y_2$, $p_5=1$, and $\Delta$ is the unbounded polygon defined by $y_1,y_2>0$, $7y_1-11y_2+2>0$ and $9y_1-9y_2+1>0$. By Proposition~\ref{prop:vertices}, $V(\Delta)$ consists of the three vertices $\{3,4\}$, $\{2,3\}$, $\{1,2\}$, at $(0,0)$, $(0,\frac19)$, $(\frac7{36},\frac{11}{36})$, and of the two vertices at infinity $\{1,5\}$ and $\{4,5\}$, in the directions $(11,7)$ and $(1,0)$. The edges of $Q$ are indexed by $1,\dots,5$ in cyclic order. The signs are given in the following table:
\[
\begin{array}{c|ccccc}
S&\{1,2\}&\{2,3\}&\{3,4\}&\{4,5\}&\{1,5\}\\ \hline
\det B_S&4&5&4&4&-5\\
\sign\det(D_S;D_{i_S})&+1&+1&+1&+1&-1
\end{array}
\]
So $\omega(S)=+1$ for all $S\in V(\Delta)$. For $e=(0,1)$, the vertices such that $e\in\cone(-B_s,\ s\in S)$ are $\{1,2\}$ and $\{3,4\}$, and for $e=(1,0)$ they are $\{1,5\}$ and $\{2,3\}$. We get $\deg(\Delta)=2$, and $\deg F=\deg(\Delta)$ for these choices of $B$ and $D$. For the heights $h=(-1,0,0,0,0)$ of Example~\ref{ex:run}, the vector $e=-B^Th=B_1=(-5,-2)$ selects $\{4,5\}$ and $\{2,3\}$. As predicted by Lemma~\ref{lem:translation}, these are the complements of the positively decorated triangles $\{1,2,3\}$ and $\{1,4,5\}$ of $T_h$. The closed polygonal line with successive vertices $B_1,\dots,B_5$ turns twice around the origin, as predicted by Proposition~\ref{prop:k2}.
\end{example}

\section{The topological degree as a half signed count}\label{sec:ends}

In this section we show that the topological degree is half a signed count, which is the key to our Descartes' rule in Section~\ref{sec:descartes} and to the bounds of Section~\ref{sec:bound}. The idea is simple. By the Newton polytope formula, $\deg F$ is the number of preimages, counted with signs, of a direction $u$ by a piecewise linear map from $\partial P$ to the sphere $S^{n-1}$. Adding the formulas for $u=e_n$ and $u=-e_n$, we get that $2\deg F$ is a signed count of the points of $\partial P$ where the first $n-1$ coordinates of this map vanish (Proposition~\ref{prop:ends-newton}). Geometrically, these points correspond to the ends of the curve $f_1=\dots=f_{n-1}=0$. The same argument, applied to the Gale dual formula, gives Proposition~\ref{prop:ends-gale}.

\begin{proposition}\label{prop:ends-newton}
Assume (H), and that neither $e_n$ nor $-e_n$ lies in a cone spanned by $n-1$ columns of $C$. Let $C'$ be the matrix formed by the first $n-1$ rows of $C$. Then, for any triangulation $T$ of $\partial P$ as in Theorem~\ref{thm:Ap},
\[
\deg F=\frac{(-1)^n}2\sum_{\substack{\tau\in T,\ \dim\tau=n-1\\ (C'_s)_{s\in\tau}\ \text{positively span}\ \R^{n-1}}}\sign\det C^\tau\cdot\sign\det A^{\tau\cup\{i_\tau\}} .
\]
By Cramer's rule, $(C'_s)_{s\in\tau}$ positively spans $\R^{n-1}$ if and only if the numbers $(-1)^r\det C'^{\tau\setminus\{s_r\}}$, $r=1,\dots,n$, are nonzero and have the same sign, where $\tau=\{s_1<\dots<s_n\}$.
\end{proposition}

This assumption holds if all the maximal minors of $C'$ are nonzero, since $\det(C^S,e_n)=\pm\det C'^S$ for $|S|=n-1$. It holds after replacing $F$ by $MF$ for a generic $M\in GL_n(\R)$, since a cone spanned by $n-1$ columns of $C$ lies in a proper subspace of $\R^n$, and a generic line through the origin meets such a cone only at the origin. This does not change the positive zeros, and it multiplies $\deg F$ by $\sign\det M$. So the vector $u$ of Theorem~\ref{thm:Ap} can be replaced by a generic change of basis. The same holds on the Gale side (Proposition~\ref{prop:ends-gale}).

\begin{proof}
Let $g=\Pi_T/|\Pi_T|\colon\partial P\to S^{n-1}$, which is defined by Lemma~\ref{lem:homotopy}, and whose degree is $\deg F$. By assumption, Lemma~\ref{lem:PL} applies both with $u=e_n$ and with $u=-e_n$: $\deg g$ is the signed count of the points of $g^{-1}(e_n)$, and also of the points of $g^{-1}(-e_n)$, where a point of an $(n-1)$-simplex $\tau$ is counted with the sign $c_\tau=(-1)^n\sign\det C^\tau\cdot\sign\det A^{\tau\cup\{i_\tau\}}$. This sign does not depend on $u$. Averaging the two counts, we get
\[
\deg F=\frac12\sum_{q\in g^{-1}(\{e_n,-e_n\})}c_{\tau(q)},
\]
where $\tau(q)$ is the simplex containing $q$.

The points of $g^{-1}(\{e_n,-e_n\})$ are the points of $\partial P$ where the first $n-1$ components of $\Pi_T$ vanish. Consider an $(n-1)$-simplex $\tau=\{s_1<\dots<s_n\}$ of $T$. A point $q=\sum_r\beta_ra_{s_r}$ of $\tau$, with $\beta_r\ge0$ and $\sum_r\beta_r=1$, is such a point if and only if $\sum_r\beta_rC'_{s_r}=0$. Assume that $q$ is such a point. Then $\Pi_T(q)=\sum_r\beta_rC_{s_r}$ is a nonzero multiple of $e_n$, by Lemma~\ref{lem:homotopy}. If some $\beta_r$ were zero, $e_n$ or $-e_n$ would lie in a cone spanned by $n-1$ columns of $C$. It follows that all the $\beta_r$ are positive. For the same reason, $\det C^\tau\neq0$ by Carath\'eodory's theorem, so the $(n-1)\times n$ matrix $C'^\tau$ has rank $n-1$, and its kernel is a line. So the vectors $C'_s$, $s\in\tau$, positively span $\R^{n-1}$, and $q$ is the only such point of $\tau$. Conversely, assume that the vectors $C'_s$, $s\in\tau$, positively span $\R^{n-1}$, and let $\beta$ be a positive vector of the kernel of $C'^\tau$, normalized so that $\sum_r\beta_r=1$. Then $q=\sum_r\beta_ra_{s_r}$ is a point of $\tau$ where the first $n-1$ components of $\Pi_T$ vanish. We get that $\tau$ contains such a point if and only if the vectors $C'_s$, $s\in\tau$, positively span $\R^{n-1}$, and in this case it contains exactly one, $q_\tau$, which lies in the relative interior of $\tau$. This gives the formula. Finally, by Cramer's rule, the kernel of $C'^\tau$ is spanned by the vector $\big((-1)^r\det C'^{\tau\setminus\{s_r\}}\big)_r$, which gives the last assertion.
\end{proof}

Geometrically, after a generic change of basis of $\R^n$, the points counted in Proposition~\ref{prop:ends-newton} correspond, through the moment map, to the ends of the curve $f_1=\dots=f_{n-1}=0$ in $\R^n_{>0}$. The positive zeros of $F$ are the zeros of $f_n$ on this curve, and their Jacobian signs alternate along each component. So $\deg F$ is half the signed count of the ends of this curve, and $|\deg F|$ is at most its number of noncompact components.

\begin{proposition}\label{prop:ends-gale}
Assume (U$_C$), (H) and $\Delta\neq\emptyset$, and that neither $e_k$ nor $-e_k$ lies in a cone spanned by $k-1$ of the vectors $-B_i$. Let $B'$ be the matrix formed by the first $k-1$ columns of $B$. Then, with the global sign of Theorem~\ref{thm:B},
\[
\deg F=\pm\frac12\sum_{\substack{S\in V(\Delta)\\ (B'_s)_{s\in S}\ \text{positively span}\ \R^{k-1}}}\omega(S).
\]
\end{proposition}

The assumption holds in particular if all the maximal minors of $B'$ are nonzero. It holds after replacing $B$ by $BM$ for a generic $M\in GL_k(\R)$, which amounts to a generic change of basis of $\ker A$. This does not change $\Delta$, and it only changes the global sign. As on the Newton side, the vertices counted in Proposition~\ref{prop:ends-gale} correspond to the ends of the curve defined in $\Delta$ by the first $k-1$ equations $\psi_1=\dots=\psi_{k-1}=0$ of the Gale system~\eqref{eq:gale-system}, and the bound of Theorem~\ref{thm:descartes}(2) below is the bound for the number of noncompact components of this curve obtained in \cite{PartI}.

\begin{proof}
By assumption, Theorem~\ref{thm:Bp} applies both with $e=e_k$ and with $e=-e_k$. By the proof of Theorem~\ref{thm:B}, the global sign $\epsilon$ only depends on $B$ and $D$, so it is the same in both cases. A vertex $S\in V(\Delta)$ contributes $\omega(S)$ to the formula for $e=\pm e_k$ if and only if $\pm e_k\in\cone(-B_s,\ s\in S)$. Averaging the two formulas, we get
\[
\deg F=\frac\epsilon2\sum_S\omega(S),
\]
where the sum runs over the vertices $S$ such that $\cone(-B_s,\ s\in S)$ meets the line $\R e_k$ outside the origin. Such a vertex contributes to exactly one of the two formulas: if $e_k$ or $-e_k$ lies in $\cone(-B_s,\ s\in S)$, then $\det B_S\neq0$ by Carath\'eodory's theorem, so this cone is simplicial and cannot contain both $e_k$ and $-e_k$.

It remains to characterize these vertices. The cone $\cone(-B_s,\ s\in S)$ meets $\R e_k$ outside the origin if and only if there is $\beta\in\R^S_{\ge0}$ such that $\sum_s\beta_sB'_s=0$ and $\sum_s\beta_sB_s\neq0$. The second condition is automatic for $\beta\neq0$. Assume that $\sum_s\beta_sB_s=0$ with $\beta\ge0$, $\beta\neq0$. Then the support $S'\subset S$ of $\beta$ is a coface, and $K_{S'}\supset K_S\neq\{0\}$, which contradicts (H) by Lemma~\ref{lem:gale-N}(2). Now, assume that $\sum_s\beta_sB'_s=0$ with $\beta\ge0$, $\beta\neq0$. Then $\sum_s\beta_s(-B_s)$ is a nonzero multiple of $e_k$. If some $\beta_s$ were zero, $e_k$ or $-e_k$ would lie in a cone spanned by $k-1$ of the vectors $-B_i$. It follows that all the $\beta_s$ are positive. For the same reason $\det B_S\neq0$, so the $k$ vectors $B'_s$, $s\in S$, span $\R^{k-1}$, and therefore they positively span $\R^{k-1}$. Conversely, if they positively span $\R^{k-1}$, they satisfy a linear relation with positive coefficients. This gives the formula.
\end{proof}

\section{A Descartes' rule of signs for the topological degree}\label{sec:descartes}

In this section we deduce from the half signed counts of Section~\ref{sec:ends} our Descartes' rule of signs for the topological degree (Theorem~\ref{thm:descartes}), which contains Theorem~\ref{thm:desc-intro} of the introduction. We first define the cyclic sign changes of a sequence and, in higher dimension, sign changes of families of vectors. We then define, on the Newton side and on the Gale side, the numbers of sign changes $V_T$ and $V_\Delta$, which play the role of the number of sign changes in the classical rule.

Let $(c_1,\dots,c_r)$ be a cyclic sequence of nonzero real numbers. Its number of \emph{cyclic sign changes} $V_{\rm cyc}(c_1,\dots,c_r)$ is the number of indices $j\in\{1,\dots,r\}$ such that $c_jc_{j+1}<0$, where $c_{r+1}=c_1$. This is the number of sign changes of the sequence $(c_1,\dots,c_r,c_1)$, and it is always even.

In higher dimension we use the following notion. We say that $m$ vectors of $\R^{m-1}$ form a \emph{sign change} if they positively span $\R^{m-1}$. For $m=2$, this means that two real numbers are nonzero and have opposite signs, so this is an ordinary sign change.

On the Newton side, let $g_j=\sum_lg_{lj}f_l$, $j=1,\dots,n-1$, be $n-1$ linear combinations of the components of $F$, and for $i\in[N]$ let $\gamma_i\in\R^{n-1}$ be the vector of the coefficients of $x^{a_i}$ in $g_1,\dots,g_{n-1}$. Let $T$ be a triangulation of $\partial P$ as in Theorem~\ref{thm:Ap}. We say that an $(n-1)$-simplex $\tau$ of $T$ is a \emph{sign change} of $(g_1,\dots,g_{n-1})$ if the $n$ vectors $\gamma_s$, $s\in\tau$, form a sign change. The number of sign changes of $(g_1,\dots,g_{n-1})$ on $T$ is the number of these simplices,
\[
V_T(g_1,\dots,g_{n-1})=\#\big\{\tau\in T:\ \dim\tau=n-1,\ \tau\ \text{is a sign change of}\ (g_1,\dots,g_{n-1})\big\},
\]
which is the sum over all the $(n-1)$-simplices $\tau$ of $T$ of $1$ if $\tau$ is a sign change and $0$ otherwise. Equivalently, $\tau$ is a sign change if and only if the truncated system $g_1|_\tau=\dots=g_{n-1}|_\tau=0$ has a positive solution, where $g_j|_\tau$ is the sum of the terms of $g_j$ with exponents in $\tau$ (see the proof of Theorem~\ref{thm:descartes}). Under (U$_A$), we can take for $T$ the boundary complex of $P$, whose $(n-1)$-simplices are the facets of $P$. By the vertex reduction, we can also take for $T$ a triangulation of $\partial P$ whose vertices are the vertices of $P$. We say that $(g_1,\dots,g_{n-1})$ is \emph{admissible} for $T$ if there is no vector $\lambda\ge0$, supported on at most $n-1$ vertices of $T$, such that
\[
\sum_i\lambda_i\gamma_i=0\qquad\text{and}\qquad\sum_i\lambda_iC_i\neq0 .
\]
This holds in particular if any $n-1$ of the vectors $\gamma_s$, $s$ a vertex of $T$, are linearly independent, which is the case for generic $g_j$ under (U$_C$). In general, admissibility holds for generic $g_j$ (see the proof of Theorem~\ref{thm:descartes}). For $n=2$, a linear combination $g$ is admissible if and only if its coefficient $g_v$ is nonzero at every vertex $v$ of $T$, since $C_v\neq0$ by (H).

On the Gale side, let $b_1,\dots,b_{k-1}\in\ker A$ be affine relations of $\A$, and for $i\in[N]$ let $\beta_i=(b_1(i),\dots,b_{k-1}(i))\in\R^{k-1}$. We say that a vertex $S\in V(\Delta)$ is a \emph{sign change} of $(b_1,\dots,b_{k-1})$ if the $k$ vectors $\beta_s$, $s\in S$, form a sign change, and we denote by
\[
V_\Delta(b_1,\dots,b_{k-1})=\#\big\{S\in V(\Delta):\ S\ \text{is a sign change of}\ (b_1,\dots,b_{k-1})\big\}
\]
the number of these vertices. We say that $(b_1,\dots,b_{k-1})$ is \emph{admissible} if there is no vector $\lambda\ge0$, supported on at most $k-1$ indices, such that $\sum_i\lambda_i\beta_i=0$ and $\sum_i\lambda_ib(i)\neq0$ for some $b\in\ker A$. This holds in particular if any $k-1$ of the vectors $\beta_i$ are linearly independent, which is the case for generic $b_j$ under (U$_A$), and admissibility holds for generic $b_j$ in general. For $k=2$, an affine relation $b$ is admissible if and only if $b(i)\neq0$ for every $i$ such that the linear form $b'\mapsto b'(i)$ on $\ker A$ is nonzero.

\begin{theorem}\label{thm:descartes}
Assume (H). Statements~1 and~3 hold, and if moreover (U$_C$) holds, statements~2 and~4 hold.
\begin{enumerate}
\item \emph{Newton side.} For every triangulation $T$ of $\partial P$ as in Theorem~\ref{thm:Ap}, and for all linear combinations $g_1,\dots,g_{n-1}$ of the components of $F$ which are admissible for $T$,
\[
|\deg F|\le\tfrac12V_T(g_1,\dots,g_{n-1}).
\]
\item \emph{Gale side.} For all admissible affine relations $b_1,\dots,b_{k-1}$ of $\A$,
\[
|\deg F|\le\tfrac12V_\Delta(b_1,\dots,b_{k-1})\quad\text{if }\Delta\neq\emptyset,\qquad\deg F=0\quad\text{if }\Delta=\emptyset.
\]
\item \emph{Two variables.} If $n=2$ and $v_1,\dots,v_r$ are the vertices of the Newton polygon in cyclic order, then for every linear combination $g=\lambda_1f_1+\lambda_2f_2$, $\lambda_1,\lambda_2\in\R$, whose coefficients $g_i=\langle(\lambda_1,\lambda_2),C_i\rangle$ at the vertices of the Newton polygon are nonzero,
\[
|\deg F|\le\tfrac12V_{\rm cyc}(g_{v_1},\dots,g_{v_r}).
\]
\item \emph{Codimension two.} If $k=2$ and $i_1,\dots,i_m$ are the edges of $\overline\Delta$ in cyclic order, then for every admissible affine relation $b$ of $\A$,
\[
|\deg F|\le\tfrac12V_{\rm cyc}(b_{i_1},\dots,b_{i_m}).
\]
\end{enumerate}
\end{theorem}

Statements~3 and~4 are the special cases $n=2$ and $k=2$ of statements~1 and~2. For $n=2$, we take for $T$ the triangulation of $\partial P$ by the edges $\{v_l,v_{l+1}\}$ of the Newton polygon. An edge is a sign change of $g$ if and only if $g_{v_l}g_{v_{l+1}}<0$, so that $V_T(g)=V_{\rm cyc}(g_{v_1},\dots,g_{v_r})$.  Similarly, for $k=2$ the vertices of $\overline\Delta$ are the pairs of consecutive edges $\{i_l,i_{l+1}\}$, and $V_\Delta(b)=V_{\rm cyc}(b_{i_1},\dots,b_{i_m})$. For $n=1$, Theorem~\ref{thm:A} gives the parity part of the classical Descartes' rule: $|\deg F|=1$ if the extreme coefficients have opposite signs, and $\deg F=0$ otherwise (Section~\ref{sec:onetwo}). In this sense statements~1 and~2 are Descartes' rules of signs for the topological degree in any dimension. The number of sign changes of a sequence of coefficients is replaced by the number of simplices of a triangulation of the boundary of the Newton polytope (the facets, in the generic case), or of vertices of $\overline\Delta$, at which the coefficient vectors form a sign change, and the bound only depends on the signs of the maximal minors of $A$ and of the coefficient matrix of $g_1,\dots,g_{n-1}$, or of $B$ and $D$.

\begin{proof}
Statements~1 and~2 are the bounds given by Propositions~\ref{prop:ends-newton} and~\ref{prop:ends-gale}, which express $\deg F$ as half a signed count. Statements~3 and~4 follow, as explained above.

1. Let $M\in GL_n(\R)$ be a matrix whose first $n-1$ rows are the vectors $(g_{1j},\dots,g_{nj})$, $j=1,\dots,n-1$. Then the first $n-1$ components of $MF$ are $g_1,\dots,g_{n-1}$, $\deg(MF)=\pm\deg F$, and $MF$ satisfies (H). The matrix $C'$ formed by the first $n-1$ rows of $MC$ has columns $\gamma_i$. A combination $\sum_i\lambda_iMC_i$ is a multiple of $e_n$ if and only if $\sum_i\lambda_iC_i$ lies in the line $L$ of the vectors $w\in\R^n$ such that $\sum_lg_{lj}w_l=0$ for $j=1,\dots,n-1$. Admissibility says exactly that $L$ meets the cones spanned by $n-1$ columns $C_s$, $s$ a vertex of $T$, only at the origin. Then neither $e_n$ nor $-e_n$ lies in such a cone spanned by columns of $MC$, which is what is used in the proof of Proposition~\ref{prop:ends-newton}, since only the simplices of $T$ are involved in Lemma~\ref{lem:PL}. For generic $g_j$, $L$ is a generic line, which meets these cones only at the origin, since each of them lies in a proper subspace of $\R^n$; thus generic $g_j$ are admissible. Proposition~\ref{prop:ends-newton} applied to $MF$ then expresses $2\deg F$ as $\pm$ a sum of $V_T(g_1,\dots,g_{n-1})$ signs, which gives the bound. Let now $\tau$ be an $(n-1)$-simplex of $T$. The truncated system $g_1|_\tau=\dots=g_{n-1}|_\tau=0$ has $n$ monomials with affinely independent exponents and coefficient vectors $\gamma_s$, $s\in\tau$. Its vector of monomials must be, up to a positive factor, the positive vector spanning the kernel of its coefficient matrix. So it has a positive solution if and only if the vectors $\gamma_s$, $s\in\tau$, positively span $\R^{n-1}$, which means that $\tau$ is a sign change.

2. If $\Delta=\emptyset$, then $\deg F=0$ by Theorem~\ref{thm:Bp}. Otherwise, we complete $b_1,\dots,b_{k-1}$ into a basis $b_1,\dots,b_k$ of $\ker A$, and take for $B$ the matrix with these columns. Its first $k-1$ columns form the matrix $B'$ with rows $\beta_i$. The rows of $B$ are the coordinates, in the dual basis, of the linear forms $b\mapsto b(i)$ on $\ker A$, and a combination of them is a multiple of $e_k$ if and only if the corresponding linear form vanishes on the hyperplane $W=\mathrm{span}(b_1,\dots,b_{k-1})$ of $\ker A$. Admissibility says exactly that the line of the linear forms vanishing on $W$ meets the cones spanned by $k-1$ of the forms $-B_i$ only at the origin. Then neither $e_k$ nor $-e_k$ lies in a cone spanned by $k-1$ of the vectors $-B_i$. For generic $b_j$, $W$ is a generic hyperplane and this line is a generic line, so generic $b_j$ are admissible. Then Proposition~\ref{prop:ends-gale} expresses $2\deg F$ as $\pm$ a sum of $V_\Delta(b_1,\dots,b_{k-1})$ signs, which gives the bound.
\end{proof}

In particular, if some linear combination of two equations in two variables has coefficients of the same sign at all the vertices of the Newton polygon, then $\deg F=0$. Note that, in statement~3, only the vertices of the Newton polygon matter, that the sign changes are counted cyclically, and that one may first take any linear combination of the equations.

\begin{example}[Example~\ref{ex:run}, continued: Descartes' rule]
For the running example, the coefficients of $f_1$ at the vertices $a_1,\dots,a_5$ are $3,-3,1,1,-3$. This sequence has $4$ cyclic sign changes, so statement~3 gives $|\deg F|\le2$, which is an equality. On the Gale side, the sum of the two columns of $B$ in Example~\ref{ex:run-gale} is the affine relation $b=(-7,9,-7,4,1)$ of $\A$, and the edges of $\overline\Delta$ are indexed by $1,\dots,5$ in cyclic order. The sequence $b$ has $4$ cyclic sign changes, and statement~4 also gives $|\deg F|\le2$.
\end{example}

\section{An upper bound for the topological degree}\label{sec:bound}

In this section we prove the fewnomial bound of Theorem~\ref{thm:C} and the description of the sets $\mathcal D(n,1)$ and $\mathcal D(n,2)$ of Theorem~\ref{thm:E}. The bound is obtained by counting the terms of the half signed counts of Section~\ref{sec:ends} with the Upper Bound Theorem. Its sharpness for $n=2$ or $k=2$ is shown by explicit constructions.

For $d\ge1$ and $N\ge d+1$, let $M(d,N)$ be as in Theorem~\ref{thm:C}. For example, we have $M(1,N)=2$, $M(2,N)=N$, $M(3,N)=2N-4$ and $M(4,N)=N(N-3)/2$. Recall that, by the Upper Bound Theorem \cite{McM}, $M(d,N)$ is the largest number of facets of a $d$-polytope with $N$ vertices and, dually, the largest number of vertices of a $d$-polytope with $N$ facets. To prove Theorem~\ref{thm:C}, we bound the number of terms of the half signed counts of Propositions~\ref{prop:ends-newton} and~\ref{prop:ends-gale}, for a well chosen change of basis.

\begin{proof}[Proof of Theorem~\ref{thm:C}]
If $n=1$, then $|\deg F|\le1$ by Section~\ref{sec:onetwo}, and if $k=1$, then $|\deg F|\le1$ by Section~\ref{sec:codim}. Assume $n\ge2$. Let $M\in GL_n(\R)$ be generic, so that Proposition~\ref{prop:ends-newton} applies to $MF$, with the boundary complex of $P$ as triangulation. By the proof of this proposition, a facet $\Gamma$ contributes a term if and only if $e_n$ or $-e_n$ lies in $\cone(MC_s,\ s\in\Gamma)$, or equivalently if the line $L=M^{-1}(\R e_n)$ meets $\cone(C_s,\ s\in\Gamma)$ outside the origin. So $2|\deg F|$ is at most the number of facets $\Gamma$ of $P$ such that $L$ meets the simplicial cone $\cone(C_s,\ s\in\Gamma)$. Since $n\ge2$, the union of a simplicial cone of dimension $n$ and of its opposite is not $\R^n$. We can therefore choose $M$ such that $L$ does not meet the cone of a given facet. Then $2|\deg F|$ is less than the number of facets of $P$, which is at most $M(n,N)$, since $P$ is an $n$-polytope with at most $N$ vertices. This gives $|\deg F|<\frac12M(n,N)$.

In the same way, assume $k\ge2$ and $\Delta\neq\emptyset$ (if $\Delta=\emptyset$, then $\deg F=0$). By Proposition~\ref{prop:ends-gale}, applied after a generic change of basis of $\ker A$, $2|\deg F|$ is at most the number of $S\in V(\Delta)$ such that a generic line meets the simplicial cone $\cone(-B_s,\ s\in S)$ outside the origin. Moreover, the line can be chosen so that it does not meet the cone of a given $S$. We get $2|\deg F|<\#V(\Delta)\le M(k,N)$, because $V(\Delta)$ is the vertex set of the simple $k$-polytope $Q$, which has at most $N$ facets.

Finally, if $k=2$ or $n=2$, then $\min(M(n,N),M(k,N))=N$, and the bound reads $|\deg F|\le\lceil N/2\rceil-1=\lfloor(N-1)/2\rfloor$.
\end{proof}

\begin{proof}[Proof of Theorem~\ref{thm:E}]
By Theorem~\ref{thm:C}, $\mathcal D(n,1)$ and $\mathcal D(1,n)$ are contained in $\{0,1\}$, and $\mathcal D(n,2)$ and $\mathcal D(2,n)$ are contained in $\{0,\dots,w\}$, where $w=\lfloor(N-1)/2\rfloor=\lfloor n/2\rfloor+1$ with $N=n+3$. The value $0$ is attained by any map with $\Delta=\emptyset$, for instance with all the columns of $C$ in an open half-space. For $n=1$, the value $1$ is attained by any polynomial whose extreme coefficients have opposite signs. For $k=1$, it is attained since $\mathcal D(n,1)=\mathcal D(1,n)$ by Theorem~\ref{thm:D}, or directly by a circuit with $B_{i_0}B_{i_1}<0$ (Section~\ref{sec:codim}).

\emph{The case $n=2$.} Let $\A$ be the set of vertices of a convex polygon with $N$ vertices, numbered counterclockwise. It satisfies (U$_A$). Let $1\le j\le w$, and let $C_r$ be a small generic perturbation of the vector $(\cos\frac{2\pi jr}N,\sin\frac{2\pi jr}N)$, $r=1,\dots,N$, so that (U$_C$) holds. Since $2\pi j/N<\pi$, the angle from $C_r$ to $C_{r+1}$ (indices modulo $N$) is close to $2\pi j/N$, hence lies in $(0,\pi)$, and the sum of these angles is $2\pi j$. So the closed polygonal line with successive vertices $C_1,\dots,C_N$ winds $j$ times around the origin, and $\deg F=j$ by Corollary~\ref{cor:n2}.

\emph{The case $k=2$.} Let $\theta_1<\dots<\theta_N$ be angles in $[0,2\pi)$, and let $D^0$ be the $N\times3$ matrix with rows $(-\cos\theta_i,-\sin\theta_i,1)$. The inequalities $\langle D^0_i,(y,1)\rangle>0$ define the polygon with $N$ sides circumscribed about the unit circle. No three of the tangent lines are concurrent, so any three rows of $D^0$ are linearly independent. Let $D=D^0G$ for a matrix $G\in GL_3(\R)$ such that $D_N=(0,0,1)$. It has the same cone $K$ up to a linear isomorphism, and $Q$ is a polygon with $N$ sides, in the cyclic order $1,\dots,N$. Let $C$ be a matrix whose rows form a basis of the orthogonal complement of the columns of $D$. Then $D$ is a Gale dual matrix of $C$, and $C$ satisfies (U$_C$). Next, for $1\le j\le w$, let $B$ be the $N\times2$ matrix whose rows are small generic perturbations of $(\cos\frac{2\pi jr}N,\sin\frac{2\pi jr}N)$, $r=1,\dots,N$, corrected so that its rows sum to zero. This is possible since the unperturbed rows sum to zero. Let $A$ be a matrix whose rows form a basis of the orthogonal complement of the columns of $B$, the first row being $(1,\dots,1)$. Its columns are of the form $(1,a_i)$ with $a_i\in\R^n$, and $B$ is a Gale dual matrix of $A$. Moreover, (U$_A$) holds since any two rows of $B$ are linearly independent. As in the case $n=2$, the closed polygonal line with successive vertices $B_1,\dots,B_N$ winds $j$ times around the origin, and $|\deg F|=j$ by Proposition~\ref{prop:k2}. The exponents are real. However, since the topological degree only depends on the oriented matroids (Theorem~\ref{cor:om}), they may be perturbed to rational numbers and then multiplied by a common integer. So integer exponents give the same set of values of the topological degree. Note that the case $k=2$ also follows from the case $n=2$ and Theorem~\ref{thm:D}.
\end{proof}

For $n,k\ge3$ we do not know whether the bound of Theorem~\ref{thm:C} is sharp. It gives $|\deg F|\le4$ for $n=k=3$ and $|\deg F|\le13$ for $n=k=4$.
\section{H-duality and symmetry}\label{sec:h}

In this section we prove the symmetry of Theorem~\ref{thm:D}. We first recall the H-duality of \cite{BDF24}, which associates with $F$ a fewnomial map $H$ with $k$ components in $k$ variables and $N$ monomials. We show that, for the map $H$, the roles of the exponents and of the coefficients of $F$ are exchanged at the level of Gale dual matrices (Lemma~\ref{lem:h-dual}). As a consequence, the Newton polytope formula for $F$ becomes the Gale dual formula for $H$, and we get $\deg H=\pm\deg F$ (Theorem~\ref{thm:h-duality}), from which Theorem~\ref{thm:D} follows.

Assume (U$_A$), (U$_C$) and $\Delta\neq\emptyset$. Let $v_0$ be an interior point of the cone $K$, so that $\langle D_i,v_0\rangle>0$ for all $i$. We put $\Lambda=\mathrm{diag}(\langle D_i,v_0\rangle)$ and $\tilde D_i=D_i/\langle D_i,v_0\rangle$. The points $\tilde D_i$ lie on the affine hyperplane $\{\langle\cdot,v_0\rangle=1\}$ of $\R^{k+1}$. Let $\delta_i\in\R^k$ be their coordinates in an affine coordinate system of this hyperplane. Following \cite[Section~5]{BDF24}, we call
\begin{equation}\label{eq:h-dual}
H=(H_1,\dots,H_k)\colon\R^k_{>0}\to\R^k,\qquad H_j(z)=\sum_{i=1}^NB_{ij}\,z^{\delta_i},
\end{equation}
an \emph{H-dual} of $F$. This is a fewnomial map with $k$ components in $k$ variables, with $N=k+n+1$ monomials, exponent set $\{\delta_1,\dots,\delta_N\}$ and coefficient matrix $B^T$. In \cite{BDF24}, the point $v_0$ is attached to a positive zero $x_0$ of $F$ by $Dv_0=(x_0^{a_i})_i$. Up to a positive factor, this means that $v_0=(y_0,1)$, where $y_0\in\Delta$ is the zero of the Gale system $\psi$ corresponding to $x_0$ (see~\eqref{eq:gale-system}). Moreover, in \cite{BDF24} the zero is $x_0=(1,\dots,1)$. This can always be achieved by the change of variables $x\mapsto x_0\star x$, $(x_0\star x)_j=x_{0,j}x_j$, and it means that the columns of $C$ sum to zero; then $Dv_0=(1,\dots,1)$. Since $\sum_iB_i=0$, every H-dual vanishes at $(1,\dots,1)$, and with this choice of $v_0$ this zero corresponds to $x_0$. It is used in \cite{BDF24} to compare the local multiplicities of $F$ at $x_0$ and of $H$ at $(1,\dots,1)$.

\begin{lemma}\label{lem:h-dual}
\begin{enumerate}
\item $\Lambda C^T$ is a Gale dual matrix of the exponent matrix of $H$, and $A^TG$ is a Gale dual matrix of its coefficient matrix $B^T$, for some $G\in GL_{n+1}(\R)$.
\item $H$ satisfies (U$_A$) and (U$_C$), and its positive chamber $\Delta_H$ is nonempty.
\item The vertices of $\overline{\Delta_H}$ are the sets $\Gamma\subset[N]$ which are facets of $\conv(\A)$.
\item The facets of the Newton polytope $\conv(\delta_1,\dots,\delta_N)$ of $H$ are the sets $S\in V(\Delta)$.
\end{enumerate}
\end{lemma}

\begin{proof}
1. The row space of the matrix with columns $(1,\delta_i)$ consists of the vectors $(\langle\tilde D_i,v\rangle)_i$, $v\in\R^{k+1}$. This is $\Lambda^{-1}\ker C$. Its orthogonal complement $\Lambda(\ker C)^\perp$ is spanned by the columns of $\Lambda C^T$. On the other hand, the kernel of $B^T$ is $(\ker A)^\perp$, which is the row space of $A$, spanned by the columns of $A^T$.

2. By~1, (U$_A$) for $H$ means that any $n$ rows of $\Lambda C^T$ are linearly independent, which is (U$_C$) for $F$. Likewise, (U$_C$) for $H$ is (U$_A$) for $F$. As in Section~\ref{sec:gale}, $\Delta_H\neq\emptyset$ if and only if the kernel of the coefficient matrix $B^T$ of $H$ contains a positive vector. This kernel is the row space of $A$, which contains $(1,\dots,1)$.

3. By Proposition~\ref{prop:vertices}, a set $\Gamma$ with $|\Gamma|=n$ is a vertex of $\overline{\Delta_H}$ if and only if the signs of the determinants of the rows of $A^TG$ indexed by $\Gamma$, followed by the row $i$, do not depend on $i\notin\Gamma$. These signs are $\sign\det G\cdot\sign\det A^{\Gamma\cup\{i\}}$. They do not depend on $i$ exactly when all the other points of $\A$ lie on the same side of the affine hull of $\Gamma$, which means that $\Gamma$ is a facet.

4. Recall that the $\delta_i$ are the coordinates of the points $\tilde D_i$ in an affine coordinate system of the hyperplane $\{\langle\cdot,v_0\rangle=1\}$. Homogenizing this coordinate system, we get an invertible $(k+1)\times(k+1)$ matrix $M'$ such that $(1,\delta_i)=M'\tilde D_i=M'D_i/\langle D_i,v_0\rangle$ for all $i$. It follows that, for $|S|=k$ and $i\notin S$,
\[
\det\big((1,\delta_s)_{s\in S},(1,\delta_i)\big)=\det M'\cdot\frac{\det(D_S;D_i)}{\prod_{s\in S}\langle D_s,v_0\rangle\cdot\langle D_i,v_0\rangle},
\]
with positive denominators. A set $S$ of size $k$ is a facet of $\conv(\delta_1,\dots,\delta_N)$ if and only if the sign of the left-hand side does not depend on $i\notin S$. This means that the signs of $\det(D_S;D_i)$, $i\notin S$, are equal, and by Proposition~\ref{prop:vertices} this means that $S\in V(\Delta)$.
\end{proof}

Items~3 and~4 express the involutivity of the H-duality. In the same way, one checks that an H-dual of $H$ is the map $F$, with its exponents replaced by their images under a projective transformation which is defined on $\conv(\A)$, and with its coefficient columns $C_i$ multiplied by the positive numbers $\langle D_i,v_0\rangle$. Neither change affects the signs of the maximal minors, so this H-dual has the same topological degree as $F$.

\begin{theorem}\label{thm:h-duality}
Assume (U$_A$), (U$_C$) and $\Delta\neq\emptyset$. Then $\deg(\Delta_H)$ is, up to sign, the sum of Theorem~\ref{thm:A} for $F$, and consequently
\[
\deg H=\pm\deg F .
\]
\end{theorem}

\begin{proof}
We apply the formula of Theorem~\ref{thm:B} to $H$, with the Gale dual matrices of Lemma~\ref{lem:h-dual}(1) and a generic vector $e_H\in\R^n$. By Lemma~\ref{lem:h-dual}(3), the sum runs over the facets $\Gamma$ of $\conv(\A)$. Since $\Lambda$ is positive, the rows of $\Lambda C^T$ indexed by $\Gamma$ have a determinant of sign $\sign\det C^\Gamma$. The second sign in $\omega(\Gamma)$ is $\sign\det G\cdot\varepsilon_\Gamma$. Finally, $e_H\in\cone(-\Lambda_{ss}C_s,\ s\in\Gamma)$ if and only if $-e_H\in\cone(C_s,\ s\in\Gamma)$. So $\deg(\Delta_H)$ is $\sign\det G$ times the sum of Theorem~\ref{thm:A} for $F$, computed with $u=-e_H$. By Theorems~\ref{thm:A} and~\ref{thm:B}, we conclude that $\deg H=\pm\deg(\Delta_H)=\pm\deg F$.
\end{proof}

The H-duality thus exchanges the two formulas: the formula of Theorem~\ref{thm:B} for $H$, over the vertices of $\overline{\Delta_H}$, is the formula of Theorem~\ref{thm:A} for $F$, over the facets of $\conv(\A)$. Note that only the signed counts are preserved. For instance, consider $x^3-2x+1=(x-1)(x^2+x-1)$, with the exponents ordered as $a_1=0$, $a_2=1$, $a_3=3$, so that $B$ has rows $2,-3,1$. It has two positive roots of opposite signs. With $v_0=(1,1)$ and $D$ with rows $(2,1),(1,1),(0,1)$, its H-dual is $2z^{2/3}-3z^{1/2}+1$, which, with $t=z^{1/6}$, is $(t-1)(2t^3-t^2-t-1)$. It also has two positive roots of opposite signs, and both topological degrees vanish. The same translation also exchanges the Viro formula for $H$ and the Newton polytope formula for $F$. To see this, multiply the coefficient $B_i$ of $z^{\delta_i}$ in~\eqref{eq:h-dual} by $t^{g_i}$, where $g\in\R^N$ is such that $u=C\Lambda g$ is generic as in Theorem~\ref{thm:A}. By Lemma~\ref{lem:translation} applied to $H$, with the Gale dual matrices of Lemma~\ref{lem:h-dual}(1), the simplices of the regular triangulation of $\{\delta_1,\dots,\delta_N\}$ induced by $g$ which are positively decorated by $B^T$ are the complements of the facets $\Gamma$ of $\conv(\A)$ such that $u\in\cone(C_s,\ s\in\Gamma)$. It follows that the sum of Theorem~\ref{thm:A} for $F$ is the Viro formula for $H$, which gives another explanation of the interpretation of Theorem~\ref{thm:A} as a Viro formula in \cite[Remark~3.9]{Degree1}.

\begin{example}[Example~\ref{ex:run}, continued: the H-dual]
We take the matrices $B$ and $D$ of Example~\ref{ex:run-gale}, and $v_0=(\frac16,\frac16,1)$, which lies in the interior of $K$ since $Dv_0=(\frac43,1,1,1,1)$. As affine coordinates on the hyperplane $\{\langle\cdot,v_0\rangle=1\}$, we take $\frac43$ times the first two coordinates. Then $\delta_1=(7,-11)$, $\delta_2=(12,-12)$, $\delta_3=(8,0)$, $\delta_4=(0,8)$, $\delta_5=(0,0)$, and
\[
H(z)=\big(-5z_1^7z_2^{-11}+7z_1^{12}z_2^{-12}-6z_1^8+4z_2^8,\ -2z_1^7z_2^{-11}+2z_1^{12}z_2^{-12}-z_1^8+1\big).
\]
Its Newton polygon is the pentagon with vertices $\delta_5,\delta_1,\delta_2,\delta_3,\delta_4$ in counterclockwise order, and its edges are the five vertices of $\overline\Delta$ (Lemma~\ref{lem:h-dual}(4)). Its coefficient columns are the rows $B_i$, and the closed polygonal line with successive vertices $B_5,B_1,B_2,B_3,B_4$ turns twice around the origin. We get $\deg H=2=\deg F$ by Corollary~\ref{cor:n2}. Note that the computation of $\deg H$ by Theorem~\ref{thm:A} is exactly the computation of $\deg F$ by Theorem~\ref{thm:B} in Example~\ref{ex:run-gale}.
\end{example}

\begin{proof}[Proof of Theorem~\ref{thm:D}]
A map with $\Delta=\emptyset$ has topological degree $0$, and such maps exist for all $(n,k)$, so $0$ belongs to both sets. Let $F$ be of type $(n,k)$, satisfying (U$_A$) and (U$_C$), with $\Delta\neq\emptyset$. Its H-dual $H$ has type $(k,n)$, it satisfies (U$_A$) and (U$_C$) by Lemma~\ref{lem:h-dual}(2), and $|\deg H|=|\deg F|$ by Theorem~\ref{thm:h-duality}. So $\mathcal D(n,k)\subset\mathcal D(k,n)$, and exchanging $n$ and $k$ gives the equality. The exponents of $H$ are real. However, since the topological degree only depends on signs of minors (Theorem~\ref{cor:om}), they may be perturbed to rational numbers and then, by the change of variables $z\mapsto z^q$, to integers. It follows that the statement also holds for maps with integer exponents.
\end{proof}

\begin{remark}
In terms of oriented matroids, the H-duality is the map $(M_A,M_C)\mapsto(M_C^*,M_A^*)$, where $M_A$ and $M_C$ are the oriented matroids of the columns of $A$ and $C$, and $M^*$ denotes the dual oriented matroid of $M$ \cite[Section~3.4]{BLSWZ}. Recall that if $M$ is realized by the columns of a matrix, then $M^*$ is realized by the rows of a Gale dual matrix. So $M_A^*$ and $M_C^*$ are the oriented matroids of the rows of $B$ and of $D$. Recall also that an oriented matroid realized by vectors $v_1,\dots,v_N$ is \emph{acyclic} if there is no linear dependence $\sum_i\lambda_iv_i=0$ with $\lambda\ge0$ and $\lambda\neq0$.

The pairs $(M_A,M_C)$ that we consider satisfy two conditions. First, $M_A$ is acyclic, since the first coordinate of the columns $(1,a_i)$ of $A$ is equal to $1$. Second, $M_C^*$ is acyclic, since this means that $\ker C$ contains a positive vector, that is, $\Delta\neq\emptyset$. The H-duality exchanges these two conditions: the first component $M_C^*$ of the image is acyclic by the second condition, and the dual $(M_A^*)^*=M_A$ of its second component is acyclic by the first condition. Moreover, applying the map twice gives back $(M_A,M_C)$, since $M^{**}=M$. So the H-duality is an involution of the set of pairs of uniform oriented matroids satisfying these two conditions. Here uniformity corresponds to (U$_A$) and (U$_C$), and it is preserved by duality. Since $\deg F$ only depends on $M_A$ and $M_C$ (Theorem~\ref{cor:om}), Theorem~\ref{thm:h-duality} says that this involution preserves $|\deg F|$.
\end{remark}

% Bibliography embedded (generated by BibTeX with amsplain)
\providecommand{\bysame}{\leavevmode\hbox to3em{\hrulefill}\thinspace}
\providecommand{\MR}{\relax\ifhmode\unskip\space\fi MR }
\providecommand{\MRhref}[2]{%
  \href{http://www.ams.org/mathscinet-getitem?mr=#1}{#2}
}
\providecommand{\href}[2]{#2}

\end{document}